\documentclass[letterpaper,12pt]{article}
\usepackage{amssymb}
\usepackage{amsmath}
\usepackage{makeidx}
\usepackage{amscd}

\newtheorem{theo}{Theorem}[section]
\newtheorem{prop}[theo]{Proposition}
\newtheorem{lem}[theo]{Lemma}
\newtheorem{cor}[theo]{Corollary}
\newtheorem{defi}[theo]{Definition}

\def \Br {{\rm{Br}}}

\def \si {{\sigma}}
\def \Ga {{\Gamma}}

\def \Pic {{\rm {Pic}}}

\def \Gal {{\rm{Gal}}}
\def \Ker {{\rm{Ker}}}
\def \Im {{\rm {Im}}}

\def \P{{\mathbb P}}

\def \dim {{\rm{dim}}}
\def \Hom {{\rm {Hom}}}
\def \End {{\rm {End}}}
\def \Pic {{\rm {Pic}}}
\def \GL {{\rm {GL}}}

\def \Aut{{\rm Aut}}

\def\ov{\overline}

\def \Z {{\mathbb Z}}
\def \Q {{\mathbb Q}}
\def \F {{\mathbb F}}

\def \Id {{\rm Id}}

\def \Tr {{\rm{Tr}}}
\def \Mat {{\rm{Mat}}}

\def \Id {{\rm{Id}}}
\def \val {{\rm{val}}}

\def\G{{\bf G}}

\def\T{{\cal T}}

\def\lra{\longrightarrow}

\def\H{{\rm H}}

\def\Tr{{\rm Tr}}

\def\O{{\cal O}}

\def\Kum{{\rm Kum}}

\def\NS{{\rm NS\,}}

\def\O{{\cal O}}

\def\sE{{\cal E}}
\def\sC{{\cal C}}

\def\si{\sigma}

\def\res{{\rm res}}
\def\val{{\rm val}}

\def\Ga{\Gamma}

\def\et{\rm{\acute et}}

\def\p{{\mathfrak p}}

\def\fA{{\mathfrak A}}

\newcommand{\bthe}{\begin{theo}}
\newcommand{\ble}{\begin{lem}}
\newcommand{\bpr}{\begin{prop}}
\newcommand{\bco}{\begin{cor}}
\newcommand{\bde}{\begin{defi}}
\newcommand{\ethe}{\end{theo}}
\newcommand{\ele}{\end{lem}}
\newcommand{\epr}{\end{prop}}
\newcommand{\eco}{\end{cor}}
\newcommand{\ede}{\end{defi}}

\title{The Brauer group of Kummer surfaces and torsion of elliptic curves}
\author{Alexei N. Skorobogatov and Yuri G. Zarhin}
\date{}
\begin{document}
\baselineskip=15pt
\maketitle

\section*{Introduction}

In this paper we are interested in computing the Brauer
group of K3 surfaces. To an element of the Brauer--Grothendieck group $\Br(X)$
of a smooth projective variety $X$ over a number field $k$
class field theory associates the corresponding Brauer--Manin
obstruction, which is a closed condition satisfied by 
$k$-points inside the topological space of adelic points of $X$,
see \cite[Ch. 5.2]{Sk}. If such a condition
is non-trivial, $X$ is a counterexample to weak approximation, and
if no adelic point satisfies this condition, $X$ is a counterexample
to the Hasse principle. The computation of $\Br(X)$ is thus a first
step in the computation of the Brauer--Manin obstruction on $X$.

Let $k$ be an arbitrary field with a separable closure $\ov k$,
$\Ga=\Gal(\ov k/k)$.
Recall that for a variety $X$ over $k$
the subgroup $\Br_0(X)\subset\Br(X)$ denotes the image of $\Br(k)$ in
$\Br(X)$, and $\Br_1(X)\subset\Br(X)$ denotes the kernel of the natural
map $\Br(X)\to \Br(\ov X)$, where $\ov X=X\times_k\ov k$.
In \cite{SZ} we
showed that if $X$ is a K3 surface over a field $k$ finitely
generated over $\Q$, then $\Br(X)/\Br_0(X)$ is finite. No general
approach to the computation of $\Br(X)/\Br_0(X)$ seems to be known;
in fact until recently there was not a single K3 surface over a number field
for which $\Br(X)/\Br_0(X)$ was known. One of the aims of this
paper is to give examples of K3 surfaces $X$ over $\Q$ such that
$\Br(X)=\Br(\Q)$.

We study a particular kind of K3 surfaces, namely Kummer
surfaces $X=\Kum(A)$ constructed from abelian surfaces $A$. Let
$\Br(X)_n$ denote the $n$-torsion subgroup of $\Br(X)$. 
Section \ref{s0} is devoted to the geometry of Kummer surfaces.
We show that there is a natural isomorphism
of $\Ga$-modules 
$\Br(\ov X)\tilde\lra\Br(\ov A)$ (Proposition \ref{2.5}). When $A$
is a product of two elliptic curves, the algebraic Brauer group
$\Br_1(X)$ often coincides with $\Br(k)$, see Proposition
\ref{2.1}.

Section \ref{s2} starts with a general remark on the \'etale
cohomology of abelian varieties that may be of independent interest
(Proposition \ref{h1}). It implies that if $n$ is an odd integer,
then for any
abelian variety $A$ the group $\Br(A)_n/\Br_1(A)_n$ is canonically
isomorphic to the quotient of $\H^2_{\et}(\ov A,\mu_n)^\Ga$ by
$(\NS(\ov A)/n)^\Ga$, where $\NS(\ov A)$ is the N\'eron--Severi
group (Corollary \ref{h2}). For any $n\geq 1$ we prove that $\Br(X)_n/\Br_1(X)_n$
is a subgroup of $\Br(A)_n/\Br_1(A)_n$, and this inclusion is
an equality for odd $n$, see Theorem \ref{h6}. We deduce that
the subgroups of elements of odd order of the transcendental Brauer groups
$\Br(X)/\Br_1(X)$ and $\Br(A)/\Br_1(A)$ are naturally isomorphic.

More precise results are obtained in Section \ref{s3} in the
case when $A=E\times E'$ is a product of two elliptic curves. In
this case for any $n\geq 1$ we have 
$$\Br(A)_n/\Br_1(A)_n=\Hom_\Ga(E_n,E'_n)/\big(\Hom(\ov E,\ov {E'})/n\big)^\Ga$$
(Proposition \ref{n2}). This gives a convenient formula for
$\Br(X)_n/\Br_1(X)_n$ when $n$ is odd. See Proposition \ref{res1}
for the case $n=2$.

In Section \ref{s4} we find many pairs of elliptic curves $E,\,E'$ over
$\Q$ such that for $A=E\times E'$ the group $\Br(A)/\Br_1(A)$ is
zero or a finite abelian 2-group. For example,
if $E$ is an elliptic curve over $\Q$ such
that for all primes $\ell$ the representation $\Ga\to
\Aut(E_\ell)\simeq \GL(2,\F_\ell)$ is surjective, then for
$A=E\times E$ we have $\Br(A)=\Br_1(A)$, whereas 
$\Br(\ov A)^\Ga\simeq\Z/2^m$ for some $m\geq 1$ (Proposition \ref{pr1}).
This shows, in particular, that the Hochschild--Serre spectral sequence
$\H^p(k,\H^q_{\et}(\ov A,\G_m))\Rightarrow \H^{p+q}_{\et}(A,\G_m))$
does not degenerate. For this $A$ the
corresponding Kummer surface $X=\Kum(A)$ has trivial Brauer group 
$\Br(X)=\Br(\Q)$ (whereas 
$\Br(\ov X)^\Ga\simeq\Z/2^m$ for some $m\geq 1$). Note that by a
theorem of W. Duke \cite{Duke} most elliptic curves over $\Q$ have
this property, see the remark after Proposition \ref{pr1}.

In Section \ref{s5} we discuss the resulting infinitely many
Kummer surfaces $X$ over $\Q$ such that $\Br(X)=\Br(\Q)$,
see (\ref{ex}-\ref{ex5}) and Examples 3 and 4. In fact most
Kummer surfaces $\Kum(E\times E')$ over $\Q$ have trivial Brauer group, see
Example A2 in Section \ref{s4}.
We also exhibit Kummer surfaces $X$ with an element of prime order
$\ell\leq 13$ in $\Br(X)$ which is not in $\Br_1(X)$. Finally, we
discuss the Brauer group of Kummer surfaces that do not necessarily
have rational points.

In the follow up paper \cite{ISZ} with Evis Ieronymou we give
an upper bound on the size of $\Br(X)/\Br_0(X)$, where $X$ is
a smooth diagonal quartic surface in $\P^3_\Q$, and give examples
when $\Br(X)=\Br(\Q)$.
The importance of K3 surfaces over $\Q$ such that
$\Br(X)=\Br(\Q)$ is that there is no Brauer--Manin obstruction, and
so the Hasse principle and weak approximation for
$\Q$-points can be tested by numerical experiments.
It would be even more interesting to get theoretical evidence
for or against the Hasse principle and weak approximation on
such surfaces.

\section{Picard and Brauer groups of Kummer surfaces
over an algebraically closed field} \label{s0}

Let $k$ be a field of characteristic zero
with an algebraic closure $\ov k$, and the absolute Galois
group $\Ga=\Gal(\ov k/k)$. Let $X$ be a
smooth and geometrically integral variety over $k$,
and let $\ov X=X\times_k\ov k$.
Let $\Br(X)=\H^2_{\et}(X,\G_m)$ be the Brauer group of $X$, and let
$\Br(\ov X)=\H^2_{\et}(\ov X,\G_m)$ be the Brauer group of $\ov X$.
For any prime number $\ell$ the Kummer sequence
$$1\to\mu_{\ell^n}\to\G_m\to\G_m\to 1$$
gives rise to the exact sequence of abelian groups
$$0\to \Pic (X)\otimes\Z/{\ell^n}\to\H_{\et}^2(X,\mu_{\ell^n})
\to \Br (X)_{\ell^n}\to 0,$$ and an exact sequence of $\Ga$-modules
\begin{equation}
0\to \Pic (\ov X)\otimes\Z/{\ell^n}\to\H_{\et}^2(\ov X,\mu_{\ell^n})
\to \Br (\ov X)_{\ell^n}\to 0. \label{ee1}
\end{equation}
If $X$ is projective, then 
$\Pic(\ov X)\otimes \Z/\ell^n=\NS(\ov X)\otimes \Z/\ell^n$, 
where $\NS(\ov X)$ is the N\'eron--Severi group of $\ov X$. So in
this case we have an exact sequence of $\Ga$-modules
\begin{equation}
0\to \NS (\ov X)\otimes\Z/{\ell^n}\to\H_{\et}^2(\ov X,\mu_{\ell^n})
\to \Br (\ov X)_{\ell^n}\to 0. \label{e1}
\end{equation}
Passing to the projective limit in (\ref{e1}) we obtain an embedding
of $\Ga$-modules
$$\NS(\ov X)\otimes\Z_\ell \hookrightarrow \H_{\et}^2(\ov X,\Z_\ell(1)).$$
The N\'eron--Severi group of an abelian variety
or a K3 surface is torsion free, so in these cases $\NS(\ov X)$
is a submodule of $\H_{\et}^2(\ov X,\Z_\ell(1))$.

\medskip

\noindent{\bf Remark.}
Let $\rho$ be the rank of $\NS(\ov{X})$, and let $b_2$ be
the second Betti number of $\ov{X}$. It is known
(\cite{Gr}, Cor. 3.4, p. 82;
\cite{EC}, Ch. 5, Remark 3.29, pp. 216--217) 
that the $\ell$-primary component $\Br(\ov{X})(\ell)\subset
\Br(\ov{X})$ is an extension of 
$\H^3_{\et}(\ov X,\Z_\ell(1))_{\rm tors}$ by
$(\Q_{\ell}/\Z_{\ell})^{b_2-\rho}$.
By Poincar\'e duality, if $X$ is a surface such that
$\NS(\ov{X})$ has no $\ell$-torsion, then
$\Br(\ov{X})(\ell) \simeq (\Q_{\ell}/\Z_{\ell})^{b_2-\rho}$. It follows
that if $X$ is an abelian variety or a K3 surface 
we have $\Br(\ov{X})\simeq (\Q/\Z)^{b_2-\rho}$.

\medskip

We write $k[X]$ for the $k$-algebra of regular functions on $X$, and
$k[X]^*$ for the group of invertible regular functions.
We state the following well known fact for future reference.

\ble \label{st}
Let $X$ be a smooth and geometrically integral variety over $k$,
and let $U\subset X$ be an open subset whose complement in $X$ has codimension
at least $2$. Then the natural restriction maps
$$k[X]\to k[U], \quad \Pic(X)\to\Pic(U), \quad \Br(X)\to \Br(U)$$
are isomorphisms.
\ele
{\em Proof} The first two statements are clear, and the last one
follows from Grothendieck's purity theorem, see 
\cite{Gr}, Cor. 6.2, p. 136. QED

\medskip

For an abelian variety $A$ we denote by
$A_n$ the kernel of the multiplication by $n$ map $[n]:A\to A$.
Let $\iota$ be the antipodal involution on $A$, $\iota(x)=-x$.
The set of fixed points of $\iota$ is $A_2$.

Assume now that $A$ is an abelian surface.
Let $A_0=A\setminus A_2$ be the complement to $A_2$, and let
$X_0=A_0/\iota$. The surface $X_0$
is smooth and the morphism $A_0\to X_0$ is a torsor under $\Z/2$. Let
$X$ be the surface obtained by blowing-up the singular points of
$A/\iota$. Then $X$ can be viewed as a smooth compactification
of $X_0$; the complement to $X_0$ in $X$ is a closed subvariety of
dimension $1$ which splits over $\ov k$ into a disjoint union of
$16$ smooth rational curves with self-intersection $-2$. We shall
call $X$ the Kummer surface attached to $A$, and write $X=\Kum(A)$.

Let $A'$ be the surface obtained by blowing-up the subscheme $A_2$ in
$A$, and let $\sigma:A'\to A$ be the resulting birational morphism. Let
$\pi:A'\to X$ be the natural finite morphism of degree 2
ramified at $X\setminus X_0$ (cf. \cite{N}). 
The set $A_2(\ov k)$ is the disjoint union of $\Ga$-orbits
$\Lambda_1, \ldots, \Lambda_r$. One may view each $\Lambda_i$ as a closed
point of $A$ with residue field $K_i$. Then $M_i=\si^{-1}(\Lambda_i)$ 
in $A'$ is the projective line $\P^1_{K_i}$ (cf. \cite{ManinCubic},
Ch. III, Thm. 2.4 and Remark 2.5). It follows that $L_i=\pi(M_i)$ is
also isomorphic to $\P^1_{K_i}$.

Since $\sigma:A'\to A$ is a monoidal transformation with smooth centre,
the induced maps 
$\si^*:\Br(A)\to\Br(A')$ and $\si^*:\Br(\ov A)\to\Br(\ov A')$
are isomorphisms by \cite{Gr}, Cor. 7.2, p. 138.
Let $Y\subset A'$ be an open subset containing $A_0$.
The composition of injective maps
$$\Br(A)\tilde\lra\Br(A')\to \Br(Y)\to\Br(A_0)$$
is an isomorphism by Lemma \ref{st}, and the same is true
after the base change from $k$ to $\ov k$. It follows that the following
restriction maps are isomorphisms:
\begin{equation}
\Br(A')\tilde\lra \Br(Y), \quad \Br(\ov A')\tilde\lra \Br(\ov Y).
\label{A'Y}
\end{equation}
This easily implies that the natural homomorphisms
$\Br_1(A)\to\Br_1(A')\to \Br_1(Y)$
are isomorphisms. We also obtain isomorphisms
\begin{equation}
\Br(A)_n/\Br_1(A)_n \tilde\lra\Br(A')_n/\Br_1(A')_n \tilde\lra
\Br(Y)_n/\Br_1(Y)_n.\label{isoo}
\end{equation}
Throughout the paper, we will freely use these
isomorphisms, identifying the corresponding groups.

\bpr \label{h5} Let $X_1\subset X$ be the complement to
the union of some of the irreducible
components of $X\setminus X_0$ (that is, some of the lines $L_i$). 
Then there is an exact sequence
\begin{equation}
0\to\Br(X)\to\Br(X_1)\to \oplus_i\ K_i^*/K_i^{*2}, \label{eh} 
\end{equation} 
where the sum is over $i$ such that $L_i\subset X\setminus X_1$.
In particular,
the restriction map $\Br(X)\to\Br(X_1)$ induces an isomorphism of the subgroups
of elements of odd order. 
The restriction map $\Br(\ov X)\tilde\lra \Br(\ov X_1)$
is an isomorphism of $\Ga$-modules. 
\epr 
{\em Proof} 
Let $Y=\pi^{-1}(X_1)$. From Grothendieck's exact sequence
(\cite{Gr}, Cor. 6.2, p. 137) we obtain the following commutative diagram
with exact rows:
\begin{equation}
\begin{array}{ccccccc}
0 &\to &\Br(X)& \to &\Br(X_1)& \to &
\oplus_i\,\H^1(L_i, \Q/\Z)\\
&&\downarrow&&\downarrow&&\downarrow\\
0 &\to &\Br(A')& \tilde\lra &\Br(Y)& \to &
\oplus_i\,\H^1(M_i, \Q/\Z)
\end{array} \label{brr}
\end{equation}
where both sums are over $i$ such that $L_i\subset X\setminus X_1$.
Recall that the restriction map  $\Br(A')\to\Br(Y)$ is an isomorphism by
(\ref{A'Y}), hence the right bottom arrow is zero. Let 
$\res_{L_i}:\Br(X_1)\to \H^1(L_i,\Q/\Z)$ 
and $\res_{M_i}:\Br(Y)\to \H^1(M_i,\Q/\Z)$
be the residue maps from (\ref{brr}). The double
covering $\pi:A'\to X$ is ramified at $L_i$, thus 
$\res_{M_i}(\pi^*\alpha)=2\res_{L_i}(\alpha)$. 
But this is zero, so that
$\res_{L_i}(\alpha)$ belongs to the injective image of
$\H_{\et}^1(\P^1_{K_i},\Z/2)$ in $\H_{\et}^1(\P^1_{K_i},\Q/\Z)$. Since
$\H_{\et}^1(\P^1_{\ov k},\Z/2)=0$ we deduce from the Hochschild--Serre spectral
sequence
$$\H^p(K_i,\H^q_{\et}(\P^1_{\ov k},\Z/2))\Rightarrow
\H^{p+q}_{\et}(\P^1_{K_i},\Z/2))$$ that $\H^1(\P^1_{K_i}, \Z/2)=K_i^*/K_i^{*2}$.
This establishes the exact sequence (\ref{eh}). 

The same theorem of Grothendieck \cite[Cor. 6.2]{Gr}
gives an exact sequence
$$0 \to \Br(\ov X) \to \Br(\ov X_1) \to 
\oplus_i\,\H^1(L_i\times_k\ov k, \Q/\Z)=0.$$
Since $L_i\times_k\ov k$ is a disjoint union of finitely many
copies of $\P^1_{\ov k}$, and $\H_{\et}^1(\P^1_{\ov k},\Q/\Z)=0$,
this implies the last statement
of the proposition. QED

\bpr \label{2.5}
The natural map $\pi^*:\Br(\ov X_0)\to \Br(\ov A_0)$ is an isomorphism,
so that the composed map $(\si^*)^{-1}\pi^*: \Br(\ov X)\to\Br(\ov A)$
is an isomorphism of $\Ga$-modules.
\epr
{\em Proof} Since $\pi:A_0\to X_0$ is a torsor under $\Z/2$
we have the Hochschild--Serre spectral sequence \cite[Thm. III.2.20]{EC}
$$\H^p(\Z/2,\H_{\et}^q(\ov A_0,\G_m))\Rightarrow\H^{p+q}_{\et}(\ov X_0,\G_m).$$
Let us compute a first few terms of this sequence.
By Lemma \ref{st} we have
$$\ov k[A_0]^*=\ov k^*, \quad \Pic(\ov A_0)=\Pic(\ov A), \quad
\Br(\ov A_0)=\Br(\ov A).$$
Since $k$ has characteristic 0, and $\Z/2$ acts trivially
on $\ov k^*$, the Tate cohomology group
$\hat\H^0(\Z/2,\ov k^*)$ is trivial.
We have $\H^1(\Z/2,\ov k^*)=\Z/2$. By the periodicity
of group cohomology of cyclic groups we obtain
$\H^2(\Z/2,\ov k^*)=0$.

We have an exact sequence of $\Ga$-modules
$$0\to A^t(\ov k)\to\Pic(\ov A)\to\NS(\ov A)\to 0,$$
where $A^t$ is the dual abelian surface.
The torsion-free abelian group $\NS(\ov A)$ embeds
into $\H_{\et}^2(\ov A,\Z_\ell(1))$, and since the antipodal
involution $\iota$ acts trivially on $\H_{\et}^2(\ov A,\Z_\ell(1))$,
it acts trivially on $\NS(\ov A)$, too. Thus $\H^1(\Z/2,\NS(\ov A))=0$,
so that $\H^1(\Z/2,\Pic(\ov A))$ is the image of
$\H^1(\Z/2,A^t(\ov k))$. Since $\iota$ acts on $A^t$
as multiplication by $-1$, we have
$$\H^1(\Z/2,A^t(\ov k))=A^t(\ov k)/(1-\iota)A^t(\ov k)=0.$$
Putting all this into the spectral sequence and using
Proposition \ref{h5} with $X_1=X_0$ we obtain an embedding
$\Br(\ov X)\hookrightarrow \Br(\ov A)$.
In order to prove that this is an isomorphism,
it suffices to check that the corresponding embeddings
$\Br(\ov X)_{\ell^n} \hookrightarrow \Br(\ov A)_{\ell^n}$
are isomorphisms for all primes $\ell$ and all positive integers $n$.
It is well known that
$b_2(\ov{X})=22$, $b_2(\ov{A})=6$ and $\rho(\ov X)=\rho(\ov A)+16$,
see, e.g., \cite{SP} or \cite{N}.
From this and the remark before Lemma \ref{st} is follows that the
orders of $\Br(\ov X)_{\ell^n}$ and $\Br(\ov A)_{\ell^n}$ are the same.
This finishes the proof. QED

\medskip

\noindent{\bf Remark 1} The same spectral sequence
gives an exact sequence of $\Ga$-modules
\begin{equation}
0\to\Z/2\to\Pic(\ov X_0)\to \Pic(\ov A)^\iota\to 0, \label{tt1}
\end{equation}
where $\Pic(\ov A)^\iota$ is the $\iota$-invariant subgroup of
$\Pic(\ov A)$. From (\ref{tt1}) we deduce the exact sequence
\begin{equation}
0\to\Z/2\to\Pic(\ov X_0)_{\rm tors}\to A^t_2\to 0. \label{tors}
\end{equation}
Let $\Z^{16}\subset \Pic(\ov X)=\NS(\ov X)$ be the lattice generated
by the classes of the 16 lines. Its saturation $\Pi$ in $\NS(\ov X)$
is called the Kummer lattice. In other words, $\Pi$ is
the subgroup of $\NS(\ov X)$ consisting of linear combinations
of the classes of the 16 lines
with coefficients in $\Q$. Since $\NS(\ov X)/\Z^{16}=\Pic(\ov X_0)$,
we have $\Pi/\Z^{16}=\Pic(\ov X_0)_{\rm tors}$.
It follows from (\ref{tors}) that $[\Pi:\Z^{16}]=2^5$.
Since the 16 lines are disjoint, and each of them has self-intersection $-2$,
the discriminant of $\Z^{16}$ is $2^{16}$.
Thus the discriminant of $\Pi$ is $2^6$, as was first observed in
\cite{SP}, Lemma 4 on p. 555, see also \cite{N}.

\noindent{\bf Remark 2} Considering (\ref{tt1}) modulo torsion
and taking into account that $\NS(\ov A)^\iota=\NS(\ov A)$ and
$\H^1(\Z/2,A^t(\ov k))=0$, we obtain
an isomorphism
$$\NS(\ov X)/\Pi=\Pic(\ov X_0)/\Pic(\ov X_0)_{\rm tors}\,
\tilde\lra\,\NS(\ov A).$$
In other words, we have an exact sequence of $\Ga$-modules
\begin{equation}
0\to\Pi\to \NS(\ov X)\buildrel{\si_*\pi^*}\over{\hbox to 10mm{\rightarrowfill}}
\NS(\ov A)\to 0.\label{is}
\end{equation}

\noindent{\bf Remark 3} Recall that $A_2$ acts on $X$ and $X_0$ compatibly
with its action on $A$ by translations, moreover,
the morphisms $\pi$ and $\si$ are $A_2$-equivariant.
Thus the isomorphism $(\si^*)^{-1}\pi^*: \Br(\ov X)\tilde\lra\Br(\ov A)$
is also $A_2$-equivariant. Since translations of an abelian variety act trivially
on its cohomology, the exact sequence (\ref{e1}) shows that
the induced action of $A_2$ on $\Br(\ov A)$ is trivial.
We conclude that the induced action of $A_2$ on $\Br(\ov X)$ is also trivial.

\medskip

Let us now assume that $X$ is the Kummer surface constructed from
the abelian surface $A=E\times E'$, where $E$ and $E'$ are elliptic
curves. For a divisor $D$ we write $[D]$ for the class of $D$ in the 
Picard group.

Let $C\subset\ov A$ be a curve, and
let $p:C\to E$, $p':C\to E'$ be the natural projections.
Then $p'_*p^*:\Pic^0(\ov E)\to \Pic^0(\ov {E'})$ defines
a homomorphism $\ov E\to \ov {E'}$.
This gives a well known isomorphism of Galois modules
\begin{equation} \label{n1}
\NS(\ov A)=\Z[e]\oplus \Z[e']\oplus \Hom(\ov E,\ov {E'}),
\end{equation}
where $e=E\times\{0\}$ and $e'=\{0\}\times E'$, and the $\Ga$-module
$\Hom(\ov E,\ov {E'})$ is realised inside $\NS(\ov A)$ as the
orthogonal complement to $\Z[e]\oplus \Z[e']$ with respect
to the intersection pairing.

For a curve $C\subset A$ we denote by $\si^{-1}C\subset A'$ its strict
transform in $A'$ (i.e. the Zariski closure of $C\cap A_0$ in $A'$).
In particular, if $C$ does not contain a point of order $2$ in $A$, then
$\si^{-1}C$ does not meet the corresponding  line in $A'$, and hence
$\pi(\si^{-1}C)$ does not meet the corresponding line in $X$.

We write the $\ov k$-points in $E_2$ as $\{o,1,2,3\}$ with the
convention that $o$ is the origin of the group law, and similarly
for $E'_2$.
The divisors $\{i\}\times E'$, $E\times\{j\}$, where $i\in E_2$,
$j\in E'_2$, are $\iota$-invariant, thus there are
rational curves $s_j$ and $l_i$ in $\ov X$ such that
$\pi$ induces double coverings
$$\si^{-1}(E\times \{j\})\lra s_j,
\quad \si^{-1}(\{i\}\times E')\lra l_i.$$
Let $l_{ij}$ be the line in $\ov X$ corresponding to the
2-division point $(i,j)\in A_2$.
Note that $\si_*\pi^*$ sends $[s_j]$, $[l_i]$, $[l_{ij}]$ to
$[e]$, $[e']$, $0$, respectively. Finally, let
$$f_1=2l_o+l_{oo}+l_{o1}+l_{o2}+l_{o3}, \quad
f_2=2s_o+l_{oo}+l_{1o}+l_{2o}+l_{3o}.$$ Consider the following
9-element Galois-invariant subsets of $\NS(\ov X)$:
$$\Lambda=\{[l_{ij}]\}, \quad \Sigma=\{[f_1],\ [f_2],\ [l_o],\ [l_i], \ [s_j]\},$$
where $i$ and $j$ take all values in $\{1,2,3\}$. Let $N_\Lambda$
(resp. $N_\Sigma$) be the subgroup of $\NS(\ov X)$ generated by
$\Lambda$ (resp. by $\Sigma$).

\bpr \label{2.1} Let $A=E\times E'$, where $E$ and $E'$ are elliptic
curves, and let $X=\Kum(A)$.

{\rm (i)} The set $\Lambda$ (resp. $\Sigma$) is a $\Ga$-invariant
basis of $N_\Lambda$ (resp. of $N_\Sigma$). There is an
exact sequence of $\Ga$-modules
\begin{equation}
0\to N_\Lambda\oplus N_\Sigma\to \NS(\ov X)\to \Hom(\ov E,\ov {E'})\to 0.\label{NSX}
\end{equation}

{\rm (ii)} We have $\Br_1(X)=\Br(k)$ in each of the following cases:

$E$ and $E'$ are not isogenous over $\ov k$,

$E=E'$ is an elliptic curve without complex multiplication over $\ov k$,

$E=E'$ is an elliptic curve which, over $\ov k$, has complex multiplication
by an order $\O$ of an imaginary quadratic field $K$, that is,
$\End(\ov E)=\O$, and, moreover, $\H^1(k,\O)=0$ (for
example, $K\subset k$).
\epr
{\em Proof} (i) Relations (3.8) on p. 3217 of
\cite{HS} easily imply that all 16 classes $[l_{ij}]$ are in $N_\Lambda+
N_\Sigma$. Then relation (3.9) of {\it loc. cit.} shows that $[s_o]$ also
belongs to $N_\Lambda+ N_\Sigma$. Recall that $[\Pi:\Z^{16}]=2^5$,
see Remark 1 above. It is easy to deduce that
the Kummer lattice $\Pi$ is generated by the 16 classes $[l_{ij}]$,
together with the differences $[s_i]-[s_j]$ and $[l_i]-[l_j]$ for all
possible pairs $(i,j)$. Thus $\Pi\subset N_\Lambda+ N_\Sigma$. 

A straightforward computation of the intersection pairing shows that
the lattice generated by $[f_1]$, $[f_2]$ and the classes of 16 lines
is freely generated by these elements, and so is of rank 18.
This lattice is contained in $N_\Lambda+N_\Sigma$, hence
$N_\Lambda+ N_\Sigma$ is also of rank 18 and is
freely generated by $\Lambda\cup\Sigma$, so that $N_\Lambda+
N_\Sigma=N_\Lambda\oplus N_\Sigma$. We have
$\si_*\pi^*([s_j])=[e]$ and $\si_*\pi^*([l_i])=[e']$
for any $i$ and $j$. It follows from the exact sequence (\ref{is})
that $(N_\Lambda\oplus N_\Sigma)/\Pi$ is a
$\Ga$-submodule of $\NS(\ov A)$ generated by $[e]$ and $[e']$.
We now obtain (\ref{NSX}) from (\ref{n1}).

(ii) If $\ov E$ and $\ov {E'}$ are not isogenous, then
$\Hom(\ov E,\ov {E'})=0$. If $E=E'$ is an elliptic curve without complex
multiplication, then $\Hom(\ov E,\ov {E'})=\Z$ is a trivial $\Ga$-module.
In the last case $\Hom(\ov E,\ov {E'})=\O$, so in all cases we have
$\H^1(k,\Hom(\ov E,\ov {E'}))=0$. By Shapiro's lemma
$\H^1(k,N_\Lambda)=\H^1(k,N_\Sigma)=0$, thus
$\H^1(k,\NS(\ov X))=0$ follows from the long exact sequence of
Galois cohomology attached to (\ref{NSX}).

The surface $X$ has $k$-points, for example, on $l_{oo}\simeq\P^1_k$.
This implies that the natural map $\Br(k)\to\Br(X)$ has a retraction, and
hence is injective. The same holds for the natural map
$\H^3_{\et}(k,\G_m)\to\H^3_{\et}(X,\G_m)$. Now
from the Hochschild--Serre spectral sequence
$\H^p(k,\H^q_{\et}(\ov X,\G_m))\Rightarrow \H^{p+q}_{\et}(X,\G_m))$
we obtain a split exact sequence
$$0\to\Br(k)\to\Br_1(X)\to \H^1(k,\Pic(\ov X))\to 0.$$
Since $\Pic(\ov X)=\NS(\ov X)$, this finishes the proof. QED

\section{On \'etale cohomology of abelian varieties and Kummer surfaces}
\label{s2}

We refer to \cite[Ch. 2]{Sk} for a general introduction to torsors.

Let $A$ be an abelian variety over $k$, and let $n\geq 1$. 
Let $\T$ be the $A$-torsor with structure group $A_n$ defined by the
multiplication by $n$ map $[n]:A\to A$. Let $[\T]$ be the class of
$\T$ in $\H^1_{\et}(A,A_n)$, and let $[\ov\T]$ be the image of
$[\T]$ under the natural map $\H^1_{\et}(A,A_n)\to\H^1_{\et}(\ov
A,A_n)^\Ga$. The cup-product defines a Galois-equivariant
bilinear pairing
$$\H^1_{\et}(\ov A,A_n)\ \times\ \Hom(A_n,\Z/n)\ \to\ 
\H^1_{\et}(\ov A,\Z/n).$$ Pairing with $[\ov\T]$ gives a homomorphism of
$\Ga$-modules 
$$\tau_A:\Hom(A_n,\Z/n)\to \H^1_{\et}(\ov A,\Z/n).$$
The following lemma is certainly well known, and is proved here for the 
convenience of the reader.

\ble \label{new}
$\tau_A$ is an isomorphism of $\Ga$-modules.
\ele
{\em Proof} The two groups have the same number
of elements, hence it is enough to prove the injectivity. A non-zero homomorphism
$\alpha:A_n\to\Z/n$ can be written as the composition of a surjection
$\beta:A_n\to\Z/m$ where $m|n$, $m\not=1$, followed by the
injection $\Z/m\hookrightarrow \Z/n$. The induced map 
$\H^1_{\et}(\ov A,\Z/m)\to \H^1_{\et}(\ov A,\Z/n)$ is injective, hence
if $[\ov\T]\cup \alpha=[\alpha_*\ov\T]=0$, then
the $\ov A$-torsor $\beta_*\ov\T$ under $\Z/m$ is trivial. A trivial
$\ov A$-torsor under a finite group is not connected,
whereas the push-forward $\beta_*\ov\T$ is canonically isomorphic
to the quotient $\ov A/\Ker(\beta)$, which is connected.
This contradiction shows that sending $\alpha$ to $[\alpha_*\ov\T]$
defines an injective homomorphism of abelian groups
$\Hom(A_n,\Z/n)\tilde\lra \H^1_{\et}(\ov A,\Z/n)$. The lemma is proved. QED

\bpr \label{h1} Let $A$ be an abelian variety over $k$, 
and let $m$, $n\geq 1$ and $q\geq 0$ be integers
such that $(n,q!)=1$. Then the natural group homomorphism
$$\H^q_{\et}(A,\mu_n^{\otimes m})\to 
\H^q_{\et}(\ov A,\mu_n^{\otimes m})^\Ga$$ 
has a section, and hence is surjective. 
\epr 
{\em Proof} We break the proof into three steps.

\smallskip

{\it Step} 1. Let $M$ be a
free $\Z/n$-module of rank $d$ with a basis $\{e_i\}_{i=1}^d$, and
let $M^*$ be the dual $\Z/n$-module with the dual basis 
$\{f_i\}_{i=1}^d$. For each $q\geq 1$ we have the identity map
$$\Id_{\wedge^q M}\in \End_{\Z/n}(\wedge^q M)=
\wedge^q M\otimes_{\Z/n} \wedge^q M^*, \quad 
\Id_{\wedge^q M}=
\sum (e_{i_1}\wedge\ldots \wedge e_{i_q})\otimes 
(f_{i_1}\wedge\ldots\wedge f_{i_q}),$$
where $i_1<\ldots< i_q$.
The multiplication law 
$(a\otimes b)\cdot(a'\otimes b')=(a\wedge a')\otimes(b\wedge b')$ turns 
the ring
$$\bigoplus_{q\geq 0} \wedge^q M\otimes_{\Z/n} \wedge^q M^*$$
into a commutative $\Z/n$-algebra. A straightforward calculation shows that
\begin{equation}
(\Id_M)^q=q!\, \Id_{\wedge^q M}.\label{calc}
\end{equation}

{\it Step} 2.
Recall that the cup-product defines a canonical isomorphism 
$$\wedge^q\H^1_{\et}(\ov A,\Z/n)\tilde\lra\H^q_{\et}(\ov A,\Z/n).$$ 
We have a natural homomorphism of $\Z/n$-modules
$$\H^q_{\et}(\ov A,\Z/n)\otimes_{\Z/n} \wedge^q (A_n)\to
\H^q_{\et}(\ov A,\wedge^q (A_n)),$$
and since the abelian group $\wedge^q (A_n)$ is a product of copies of $\Z/n$,
this is clearly an isomorphism.

Write $M=\H^1_{\et}(\ov A,\Z/n)$ and use 
Lemma \ref{new} to identify $A_n$ with $M^*$.
We obtain an isomorphism of $\Z/n$-modules
\begin{equation}\H^q_{\et}(\ov A,\wedge^q (A_n))=
\H^q_{\et}(\ov A,\Z/n)\otimes_{\Z/n} \wedge^q (A_n)=
\wedge^q M\otimes_{\Z/n}\wedge^q M^*.\label{coucou}
\end{equation}
The cup-product in \'etale cohomology gives rise to the map
$$\H^1_{\et}(\ov A,A_n)^{\otimes q}\to \H^q_{\et}(\ov A,A_n^{\otimes q})
\to\H^q_{\et}(\ov A,\wedge^q(A_n)),$$ and we denote by $\wedge^q[\ov\T]\in
\H^q_{\et}(\ov A,\wedge^q(A_n))^\Ga$ the image of the product of $q$
copies of $[\ov\T]$. Now (\ref{calc}) says that the isomorphism (\ref{coucou})
identifies $\wedge^q[\ov\T]$ with $q!\, \Id_{\wedge^q M}$.

We have an obvious commutative diagram of $\Ga$-equivariant pairings,
where the vertical arrows are isomorphisms:
$$
\begin{array}{ccccc}
\wedge^q M\otimes_{\Z/n}\wedge^q M^*&\times&\Hom(\wedge^qM^*,\Z/n)
&\to&\wedge^q(M)\\
\downarrow&&||&&\downarrow\\
\H^q_{\et}(\ov A,\Z/n)\otimes_{\Z/n}\wedge^q M^*&\times&\Hom(\wedge^qM^*,\Z/n)
&\to&\H^q_{\et}(\ov A,\Z/n)\\
\downarrow&&\downarrow&&||\\
\H^q_{\et}(\ov A,\wedge^q(A_n))&\times&
\Hom(\wedge^q(A_n),\Z/n)&\to &\H^q_{\et}(\ov A,\Z/n)
\end{array}
\label{pair}
$$
The pairing with the $\Ga$-invariant element
$\wedge^q[\ov\T]$ gives a homomorphism of $\Ga$-modules
$$\Hom(\wedge^q(A_n),\Z/n)\tilde\lra \H^q_{\et}(\ov A,\Z/n),$$
which is $q!$ times the canonical isomorphism
$\Hom(\wedge^qM^*,\Z/n)\tilde\lra \wedge^qM$.
By assumption $q!$ is invertible in $\Z/n$, so
this is an isomorphism of $\Ga$-modules.
Tensoring with the $\Ga$-module $\mu_n^{\otimes m}$ we obtain
an isomorphism of $\Ga$-modules
$$\Hom(\wedge^q(A_n),\mu_n^{\otimes m})\tilde\lra 
\H^q_{\et}(\ov A,\mu_n^{\otimes m}).$$
Therefore, pairing with $\wedge^q[\ov\T]$ gives rise to an isomorphism 
of abelian groups
\begin{equation}
\Hom_\Ga(\wedge^q(A_n),\mu_n^{\otimes m})\tilde\lra 
\H^q_{\et}(\ov A,\mu_n^{\otimes m})^\Ga.\label{iden}
\end{equation}

{\it Step} 3. The cup-product in \'etale cohomology gives rise to the map
$$\H^1_{\et}(A,A_n)^{\otimes q}\to \H^q_{\et}(A,A_n^{\otimes q})
\to\H^q_{\et}(A,\wedge^q(A_n)),$$ and we denote by $\wedge^q[\T]\in
\H^q_{\et}(A,\wedge^q(A_n))$ the image of the product of $q$
copies of $[\T]$.
There is a natural pairing of abelian groups
$$\H^q_{\et}(A,\wedge^q(A_n))\times\Hom_\Ga(\wedge^q(A_n),\mu_n^{\otimes
m}) \to\H^q_{\et}(A,\mu_n^{\otimes m}).$$ Pairing with
$\wedge^q[\T]$ induces a map
$\Hom_\Ga(\wedge^q(A_n),\mu_n^{\otimes m})\to
\H^q_{\et}(A,\mu_n^{\otimes m})$
such that the composition
$$\Hom_\Ga(\wedge^q(A_n),\mu_n^{\otimes m})\to
\H^q_{\et}(A,\mu_n^{\otimes m})\to \H^q_{\et}(\ov A,\mu_n^{\otimes m})^\Ga$$
is the isomorphism (\ref{iden}). This proves the proposition. QED

\medskip

In some cases the condition $(n,q!)=1$ can be dropped, 
see Corollary \ref{z2} below. 

\bco \label{h2} Let $n$ be an odd integer. Then the images of the 
groups $\Br(A)_{n}$ and
$\H^2_{\et}(\ov A,\mu_{n})^\Ga$ in $\Br(\ov A)_{n}^\Ga$ coincide,
so that we have an isomorphism 
$$\Br(A)_n/\Br_1(A)_n \simeq
\H^2_{\et}(\ov A,\mu_n)^\Ga/(\NS(\ov A)/n)^\Ga.$$
\eco
{\em Proof.} The Kummer sequences for $A$ and $\ov A$ give rise to the 
following obvious commutative diagram
with exact rows, cf. (\ref{e1}):
$$\begin{array}{ccccccccc}
0&\to &(\NS(\ov A)/n)^\Ga&\to&\H_{\et}^2(\ov A,\mu_{n})^\Ga&\to 
&\Br (\ov A)_{n}^\Ga&&\\
&&&&\uparrow\downarrow&&\uparrow&&\\
&&&&\H_{\et}^2(A,\mu_{n})&\to &\Br (A)_{n}&\to&0
\end{array}$$
The downward arrow is the section of Proposition \ref{h1}. Both
statements follow from this diagram. QED

\bthe \label{h6} Let $A$ be an abelian surface, and let $X=\Kum(A)$.
Then $\pi^*$ defines an embedding
$$ \Br(X)_n/\Br_1(X)_n\hookrightarrow \Br(A)_n/\Br_1(A)_n,$$
which is an isomorphism if $n$ is odd. The subgroups of 
elements of odd order of the transcendental Brauer groups 
$\Br(X)/\Br_1(X)$ and $\Br(A)/\Br_1(A)$
are isomorphic.
\ethe
{\em Proof} By Proposition \ref{2.5} we have the commutative diagram
\begin{equation}\begin{array}{ccc}
\Br(X)_n&\lra &\Br(A)_n\\
\downarrow&&\downarrow\\
\Br(\ov X)_n&\tilde\lra &\Br(\ov A)_n
\end{array}\label{dd1}
\end{equation}
which implies the desired embedding.
Now assume that $n$ is odd. We can write
$$\Br(A)_n=\Br(A)_n^+\oplus\Br(A)_n^-,$$
where $\Br(A)_n^+$ (resp. $\Br(A)_n^-$) is the $\iota$-invariant
(resp. $\iota$-antiinvariant) subgroup of $\Br(A)_n$. The involution
$\iota$ acts trivially on $\H^2_{\et}(\ov A,\mu_{\ell^m})$
for any $\ell$ and $m$, hence by (\ref{e1})
it also acts trivially on $\Br(\ov A)$. It follows that for odd $n$
the image of $\Br(A)_n^-$ in $\Br(\ov A)$ is zero. This gives
an isomorphism $$\Br(A)_n/\Br_1(A)_n=\Br(A)_n^+/\Br_1(A)_n^+.$$
Thm. 1.4 of \cite{ISZ} states that if $Y\to X$ is a finite flat
Galois covering of smooth geometrically irreducible varieties
with Galois group $G$, and $n$ is coprime to $|G|$, then
the natural map $\Br(X)_n\to\Br(Y)_n^G$ is an isomorphism.
We apply this to the double covering $\pi:A'\to X$.
Taking into account the isomorphism $\Br(A)=\Br(A')$ we
obtain the following commutative diagram
\begin{equation}\begin{array}{ccc}
\Br(X)_n&\tilde\lra &\Br(A)_n^+\\
\downarrow&&\downarrow\\
\Br(\ov X)_n&\tilde\lra &\Br(\ov A)_n
\end{array}\label{dd2}
\end{equation}
Our first statement follows. 
The second statement follows from the first one once we
note that an element of odd order in 
$\Br(X)/\Br_1(X)$ comes from $\Br(X)_n$ for some odd $n$.
QED

\section{The case of product of two elliptic curves} \label{s3}

We now assume that $A=E\times E'$ is the product of 
two elliptic curves. In this case we can prove the same
statement as in 
Corollary \ref{h2} but without the assumption on $n$.

The K\"unneth formula (see \cite{EC}, Cor. VI.8.13) gives a 
direct sum decomposition of $\Ga$-modules
$$\H_{\et}^2(\ov A,\Z/n)=\H_{\et}^2(\ov E,\Z/n)\oplus
\H_{\et}^2(\ov {E'},\Z/n)\oplus \H_{\et}^2(\ov A,\Z/n)_{\rm prim}\ ,$$
where 
$$\H_{\et}^2(\ov A,\Z/n)_{\rm prim}=\H_{\et}^1(\ov E,\Z/n)\otimes
\H_{\et}^1(\ov {E'},\Z/n)$$
is the {\it primitive\ } subgroup of $\H_{\et}^2(\ov A,\Z/n)$.
On twisting with $\mu_n$ we obtain the
decomposition of $\Ga$-modules
$$\H_{\et}^2(\ov A,\mu_n)=\Z/n\oplus
\Z/n\oplus \H_{\et}^2(\ov A,\mu_n)_{\rm prim}\ ,$$
where 
$$\H_{\et}^2(\ov A,\mu_n)_{\rm prim}=\H_{\et}^1(\ov E,\Z/n)\otimes
\H_{\et}^1(\ov {E'},\mu_n).$$
The canonical isomorphism
$\tau_E:\Hom(E_n,\Z/n)\tilde\lra \H^1_{\et}(\ov E,\Z/n)$
from Lemma \ref{new} gives an isomorphism of $\Ga$-modules
$$\H_{\et}^2(\ov A,\mu_n)_{\rm prim}=\Hom(E_n\otimes E'_n,\mu_n).$$ 
Using the Weil pairing we obtain a canonical isomorphism
$$\H_{\et}^1(\ov {E'},\mu_n)=\Hom(E'_n,\mu_n)=E'_n.$$
Combining all this gives canonical isomorphisms of $\Ga$-modules
\begin{equation}
\H_{\et}^2(\ov A,\mu_n)=\Z/n\oplus\Z/n\oplus
\H_{\et}^2(\ov A,\mu_n)_{\rm prim}=\Z/n\oplus\Z/n\oplus
\Hom(E_n,E'_n).
\label{e2}
\end{equation}

Let $p:A\to E$ and $p':A\to E'$ be the natural projections.
The multiplication by $n$ map $E\to E$ defines an 
$E$-torsor $\T$ with structure group $E_n$. We define
the $E'$-torsor $\T'$ similarly. The pullbacks $p^*\T$ and $p'^*\T$
are $A$-torsors with structure groups $E_n$ and $E'_n$,
respectively.
Let $[\T]\boxtimes[\T']$ be the product of $p^*[\T]$ and $p'^*[\T']$
under the pairing
$$\H^1_{\et}(A,E_n) \times \H^1_{\et}(A,E'_n) \to 
\H^2_{\et}(A,E_n\otimes E'_n).$$
Consider the natural pairing 
$$\H^2_{\et}(A,E_n\otimes E'_n)\times 
\Hom_\Ga(E_n\otimes E'_n,\mu_n)\to \H^2_{\et}(A,\mu_n).$$
Let
$$\xi:\Hom_\Ga(E_n\otimes E'_n,\mu_n)\to 
\H^2_{\et}(A,\mu_n)$$
be the map defined by pairing with $[\T]\boxtimes[\T']$.

The following map is defined by the base change from
$k$ to $\ov k$ followed by the K\"unneth projector
to the primitive subgroup:
$$\eta:\H^2_{\et}(A,\mu_n)\to \H^2_{\et}(\ov A,\mu_n)^\Ga \to 
\H^2_{\et}(\ov A,\mu_n)^\Ga_{\rm prim}=
\Hom_\Ga(E_n\otimes E'_n,\mu_n).$$

\ble \label{z1}
We have $\eta\circ\xi=\Id$. In particular, $\eta$
has a section, and hence is surjective.
\ele
{\it Proof} We must check that the composed map
\begin{equation}
\Hom(E_n\otimes E'_n,\Z/n)\to\H^2_{\et}(\ov A,\Z/n)\to 
\H^1_{\et}(\ov E,\Z/n)\otimes \H^1_{\et}(\ov {E'},\Z/n)
\label{18}\end{equation}
defined by pairing with the image of $[\T]\boxtimes[\T']$
in $\H^2_{\et}(\ov A,E_n\otimes E'_n)$ followed
by the K\"unneth projector to the primitive subgroup, 
is the isomorphism
$$\tau_{E}\otimes\tau_{E'}:\Hom(E_n\otimes E'_n,\Z/n)
\tilde\lra\H^1_{\et}(\ov E,\Z/n)\otimes \H^1_{\et}(\ov {E'},\Z/n)$$
(cf. Lemma \ref{new}). Note that the first arrow in (\ref{18}) is 
$\tau_{E}\otimes\tau_{E'}$ followed by the composed map
\begin{equation}
\H^1_{\et}(\ov E,\Z/n)\otimes \H^1_{\et}(\ov {E'},\Z/n)
\to \H^1_{\et}(\ov A,\Z/n)\otimes \H^1_{\et}(\ov A,\Z/n)
\to \H^2_{\et}(\ov A,\Z/n),\label{19}
\end{equation}
where the first arrow is $p^*\otimes p'^*$, 
and the second one is the cup-product.
By \cite[Cor. VI.8.13]{EC} the composition of (\ref{19}) with
the K\"unneth projector is the identity, hence
the composed map in (\ref{18}) is $\tau_{E}\otimes\tau_{E'}$. QED

\bco \label{z2}
For $A=E\times E'$ and any $n\geq 1$ the natural map 
$$\H^2_{\et}(A,\mu_n)\to \H^2_{\et}(\ov A,\mu_n)^\Ga$$
has a section, and hence is surjective.
\eco
{\it Proof} We have a canonical map
$p^*:\H^2_{\et}(E,\mu_n)\to \H^2_{\et}(A,\mu_n)$.
By K\"unneth decomposition and Lemma \ref{z1}
it is enough to check that
$$\H^2_{\et}(E,\mu_n)\to\H^2_{\et}(\ov E,\mu_n)^\Ga$$
has a section (and similarly for $E'$). 
The Kummer sequences for $E$ and $\ov E$
give a commutative diagram
$$\begin{array}{ccccc}
0&\to &\Z/n&\tilde\lra&\H^2_{\et}(\ov E,\mu_n)\\
&&\uparrow&&\uparrow\\
0&\to &\Pic(E)/n&\lra&\H^2_{\et}(E,\mu_n)
\end{array}$$
The left vertical arrow is given by the degree map
$\Pic(E)\to\Z$.
It has a section that sends $1\in\Z/n$ to the class
of the neutral element of $E$ in $\Pic(E)/n$. QED

\medskip

\noindent{\bf Remark}. This shows that if an abelian variety $A$
is a product of elliptic curves, then
the condition on $n$ in Proposition \ref{h1} is superfluous. 

\medskip

Recall that the natural map $\Hom(\ov E,\ov {E'})/n\to
\Hom(E_n,E'_n)$ is injective \cite[p. 124]{Milne}. Write
$\Hom(E,E')=\Hom_\Ga(\ov E,\ov {E'})$ for the group of homomorphisms
$E\to E'$.

\bpr \label{n2} For $A=E\times E'$ we have a canonical isomorphism
of  $\Ga$-modules 
$$\Br(\ov A)_n=\Hom(E_n,E'_n)/(\Hom(\ov E,\ov {E'})/n),$$
and a canonical isomorphism
of abelian groups
$$\Br(A)_n/\Br_1(A)_n=\Hom_\Ga(E_n,E'_n)/(\Hom(\ov E,\ov {E'})/n)^\Ga.$$  
\epr 
{\em Proof} 
The Kummer sequences for $A$ and $\ov A$ give rise to
the commutative diagram
$$\begin{array}{ccccccccc}
0&\to &\NS(\ov A)/n&\to&\H_{\et}^2(\ov A,\mu_{n})&\to &\Br (\ov A)_{n}&\to&0\\
&&&&\uparrow&&\uparrow&&\\
&&&&\H_{\et}^2(A,\mu_{n})&\to &\Br (A)_{n}&\to&0
\end{array}$$
Using (\ref{n1}) and (\ref{e2}) we rewrite this diagram
as follows:
$$\begin{array}{ccccccccc}
0&\to &\Hom(\ov E,\ov {E'})/n&\to&\Hom(E_n,E'_n)&\to &\Br (\ov A)_{n}&\to&0\\
&&&&\uparrow&&\uparrow&&\\
&&&&\H_{\et}^2(A,\mu_{n})&\to &\Br (A)_{n}&\to&0
\end{array}$$
The upper row here is the first isomorphism of the proposition. 
From Lemma \ref{z1} we deduce the commutative diagram
$$\begin{array}{ccccccccc}
0&\to &\big(\Hom(\ov E,\ov {E'})/n\big)^\Ga&\to&\Hom_\Ga(E_n,E'_n)&\to &\Br (\ov A)_{n}^\Ga&&\\
&&&&\uparrow\downarrow&&\uparrow&&\\
&&&&\H_{\et}^2(A,\mu_{n})&\to &\Br (A)_{n}&\to&0
\end{array}$$
where the left upward arrow is $\eta$, and the 
downward arrow is $\xi$.
The second isomorphism of the proposition is 
a consequence of the commutativity of this diagram.
QED

\bigskip

Until the end of this section we assume that
$n=2$ and the points of order 2 of $E$ and $E'$ are defined over $k$, i.e.
$E_2\subset E(k)$ and $E'_2\subset E'(k)$. 
The above considerations can then be made more
explicit. (This construction
was previously used in \cite{SS}, Appendix A.2, see also \cite{HS},
Sect. 3.2). In this case
$$\Br(\ov A)_2=\Br(\ov A)_2^\Ga=\Br(A)_2/\Br_1(A)_2.$$
Using the Weil pairing the map $\xi$ gives rise to the map
$E_2\otimes E'_2\to \Br(A)_2$ whose image maps 
surjectively onto $\Br(\ov A)_2$.
The elements of $\Br(A)_2$ obtained in this way
can be given by symbols as follows.
The curves $E$ and $E'$ can be
given by their respective equations
$$y^2=x(x-a)(x-b), \quad v^2=u(u-a')(u-b'),$$
where $a$ and $b$ are distinct non-zero elements of $k$, and
similarly for $a'$ and $b'$. 
The multiplication by $2$ torsor $E\to E$ corresponds to the 
biquadratic extension of the function field $k(E)$ 
given by the square roots of
$(x-a)(x-b)$ and $x(x-b)$, see, e.g., \cite{Knapp}, Thm. 4.2, p. 85. We choose $e_1=(0,0)$ and $e_2=(a,0)$ as a basis of $E_2$, 
and $e'_1=(0,0)$ and $e'_2=(a',0)$ as
a basis of $E'_2$; this gives rise to an obvious basis of $E_2\otimes E'_2$.
The four resulting Azumaya algebras on $A$ are written as follows:
\begin{equation}
\big((x-\mu)(x-b),(u-\nu)(u-b')\big),\quad
\mu\in\{0,\,a\}, \ \nu\in\{0,\,a'\}.\label{four}
\end{equation}
We note that the specialisation of any of these algebras at the 
neutral element of $A$ is $0\in\Br(k)$.
By the above, the classes of the algebras (\ref{four}) 
in $\Br(\ov A)$ generate $\Br(\ov A)_2$. 

The antipodal involution $\iota$ sends
$(x,y)$ to $(x,-y)$, and $(u,v)$ to $(u,-v)$, hence the Kummer surface $X=\Kum(A)$
is given by the affine equation
\begin{equation}
z^2=x(x-a)(x-b)y(y-a')(y-b').\label{e8}
\end{equation}
We denote by $A_{\mu,\nu}$ the class in $\Br(k(X))$
given by the corresponding symbol (\ref{four}).

For $A=E\times E'$ it is convenient to replace 
$X_0\subset X$ by a larger open subset. 
Let us denote by $E_2^\sharp$ the
set of $\ov k$-points of $E$ of {\it exact} order 2; in other words,
$E_2$ is the disjoint union of $\{0\}$ and $E_2^\sharp$.
Define $W\subset X$ as the complement to the 9 lines
that correspond to the points of $E_2^\sharp\times{E'}_2^\sharp$. 
The line $l_{oo}=\pi(\sigma^{-1}(0))$,
where $0\in A(k)$ is the neutral element, is contained in $W$.
Choose a $k$-point $Q$ on $l_{oo}$, and
denote by $\Br(W)^0$ the subgroup of $\Br(W)$
consisting of the elements that specialise to $0$ at $Q$.
Since $\Br(\P^1_k)=\Br(k)$, we see that $\Br(W)^0$ 
is the kernel of the restriction map $\Br(W)\to \Br(l_{oo})$,
hence $\Br(W)^0$ does not depend on the choice of $Q$.

\ble
We have $A_{\mu,\nu}\in \Br(W)^0_2$ for any 
$\mu\in\{0,\,a\}$ and $\nu\in\{0,\,a'\}$.
\ele 
{\it Proof} We have $A_{\mu,\nu}\in \Br(W)$ by \cite{SS}, Lemma A.2. 
Every $A_{\mu,\nu}$ lifts to an element of $\Br(A)$ with value $0$ 
at the neutral element of $A$, hence $A_{\mu,\nu}\in\Br(W)^0$. 
Since $A_{\mu,\nu}$ is the class of a quaternion
algebra, we have $A_{\mu,\nu}\in\Br(W)^0_2$. QED

\medskip

The map $(\mu,\nu)\mapsto A_{\mu,\nu}$ defines a group homomorphism
$\omega: E_2\otimes E'_2\to\Br(W)^0_2$.

\bpr \label{z3}
Assume one of the conditions of Proposition {\rm \ref{2.1} (ii)}.
Then we have 

{\rm (i)} $\Br_1(W)=\Br(k)$;

{\rm (ii)} $\Im(\omega)=\Br(W)^0_2$;

{\rm (iii)} $\Ker(\omega)=\Hom(\ov E,\ov {E'})/2$.
\epr
{\it Proof} (Cf. \cite{SS}, App. A2.)
(i) We have
$\ov k[W]^*=\ov k^*$, as it follows from $\ov k[A_0]^*=\ov k^*$. 
We also have $\Br_0(W)=\Br(k)$ since $W$ has a $k$-point.
Then the Hochschild--Serre spectral sequence 
$\H^p(k,\H^q_{\et}(\ov W,\G_m))\Rightarrow \H^{p+q}_{\et}(W,\G_m))$
shows that it is enough to prove 
$\H^1(k,\Pic(\ov W))=0$.
In the notation of Proposition \ref{2.1} we have 
$\Pic(\ov W)=\Pic(\ov X)/N_\Lambda$, hence
there is an exact sequence of $\Ga$-modules analogous to (\ref{NSX}):
$$0\to N_\Sigma\to \Pic(\ov W)\to \Hom(\ov E,\ov {E'})\to 0.$$
By Shapiro's lemma $\H^1(k,N_\Lambda)=0$, thus, under the assumptions
of Proposition \ref{2.1} (ii),
$\H^1(k,\Pic(\ov W))=0$ follows from the long exact sequence of
Galois cohomology.

(ii) Let $\fA\subset\Br(W)$ be the four-element set $\{A_{\mu,\nu}\}$,
and let $\ov\fA$ be the image of $\fA$ in $\Br(\ov W)$. By Proposition \ref{h5} we have $\Br(\ov W)=\Br(\ov X)$, thus we can think of $\ov\fA$
as a subset of $\Br(\ov X)$. The image of $\ov\fA$ under
the isomorphism $(\si^*)^{-1}\pi^*: \Br(\ov X)\to\Br(\ov A)$
from Proposition \ref{2.5} generates $\Br(\ov A)_2$, hence
$\ov\fA$ generates $\Br(\ov W)_2$. Therefore,
any $\alpha\in \Br(W)^0_2$ can be written as 
$$\alpha=\beta+\sum \delta_{\mu,\nu}A_{\mu,\nu},$$
where $\delta_{\mu,\nu}\in\{0,\,1\}$,
and $\beta\in\Br_1(W)$ has value zero at $Q$. It remains to apply (i).

(iii) By part (i) the natural map $\Br(W)^0 \to
\Br(\ov W)$ is injective, and we have just seen that
the latter group is naturally isomorphic to 
$\Br(\ov A)$. Now our statement follows from
the first formula of Proposition \ref{n2}.
QED

\medskip

We now calculate the residues of the $A_{\mu,\nu}$ at 
the 9 lines of $X\setminus W$ (cf. Proposition
\ref{h5} and its proof).

\ble \label{resi}
The residues of $A_{a,a'}$, $A_{a,0}$, $A_{0,a'}$, $A_{00}$,   
at the lines $l_{00}$, $l_{0,a'}$, $l_{a,0}$, $l_{a,a'}$,
written in this order, are the classes in $k^*/k^{*2}$
represented by the entries of the following matrix:
\begin{equation}\left(\begin{array}{cccc}
1 & ab & a'b' & -aa'\\
ab & 1 & aa' &a'(a'-b')\\
a'b' & aa' & 1 & a(a-b)\\
-aa'&a'(a'-b')&a(a-b)&1
\end{array}\right) \label{res}\end{equation}
For any $\mu\in\{0,\,a\}$ and $\nu\in\{0,\,a'\}$ 
the product of residues of
$A_{\mu\nu}$ at the three lines $l_{ij}$, $i\not=0$, $j\not=0$, with fixed
first or second index, is $1\in k^*/k^{*2}$. 
\ele
{\em Proof} 
We write $\res_{ij}$ for the residue at $l_{ij}$.
The local ring $O\subset k(X)$ of $l_{ij}$ is a discrete valuation ring with
valuation $\val:k(X)^*\to\Z$. For $f,\,g\in O\setminus\{0\}$ 
the residue of
$(f,g)$ at $l_{ij}$ is computed by the following rule: if $\val(f)=\val(g)=0$,
then $\res_{ij}\big((f,g)\big)$ is trivial, and if $\val(f)=0$, $\val(g)=1$,
then $\res_{ij}\big((f,g)\big)$ is the class in 
$k(l_{ij})^*/k(l_{ij})^{*2}$ of
the reduction of $f$ modulo the maximal ideal of $O$. In our case this class
will automatically be in $k^*/k^{*2}$.

Let us calculate the residues of $A_{00}=\big(x(x-b),x(y-b')\big)$.
Using the above rule we obtain
$$\res_{0,a'}(A_{00})=a'(a'-b'),\quad \res_{a,0}(A_{00})=a(a-b), \quad
\res_{a,a'}(A_{00})=1.$$
Using equation (\ref{e8}) and the relation
$(r,-r)=0$ for any $r\in k(X)^*$ we can write
$$A_{00}=\big(x(x-b),-(x-a)(y-a')\big).$$
The residue of $A_{00}$ at $l_{00}$ is then
the value of $-(x-a)(y-a')$ at $x=y=0$, that is, $-aa'$.
We thus checked the last row of (\ref{res}). 
The residue of $A_{00}$ at $l_{0,b'}$ is the class of $a(b'-a')$,
which shows that the product of residues of $A_{00}$ at $l_{00}$,
$l_{0,a'}$ and $l_{0,b'}$ is 1. The calculations in all other cases are
quite similar. QED

\medskip

\noindent{\bf Question}. Is there a conceptual explanation
of the symmetry of (\ref{res})?

\medskip

Let $r$ be the rank of $\Hom(\ov E,\ov {E'})$, and let $d$
be the dimension of the kernel of the homomorphism
$$\Hom(E_2,E'_2)=E_2\otimes E'_2\simeq(\Z/2)^4\to (k^*/k^{*2})^4$$
given by the matrix $(\ref{res})$. 

\bpr \label{res1} Let $X=\Kum(E\times E')$, where $E$ and $E'$ 
are elliptic curves with rational $2$-torsion points. 
Assume one of the conditions of Proposition {\rm \ref{2.1} (ii)}.
Then 
$$\dim_{\F_2}\ \Br(X)_2/\Br(k)_2=d-r.$$ In
particular, if $E=E'$ and 
$d=1$,
then $\Br(X)_2=\Br(k)_2$. 
\epr 
{\em Proof} This follows from Proposition \ref{z3} (ii) and (iii), and
Lemma \ref{resi}. QED

\medskip

Note that if $E=E'$ and $d=1$, then $\End(\ov E)$ has rank 1, so 
that $E$ has no complex multiplication over $\ov k$.

\section{Brauer groups of abelian surfaces} \label{s4}

In the rest of this paper we discuss abelian surfaces of the following
types:

\medskip

\noindent{\bf (A)} $A=E\times E'$, where the elliptic curves $E$ and
$E'$ are not isogenous over $\ov k$.

\noindent{\bf (B)} $A=E\times E$, where $E$ has no complex
multiplication over $\ov k$.

\noindent{\bf (C)} $A=E\times E$, where $E$ has complex
multiplication over $\ov k$.

\medskip

{\bf Case A.} In case A the N\'eron--Severi group $\NS(\ov A)$ is
freely generated by the classes $E\times\{0\}$ and $\{0\}\times E'$,
hence $\H_{\et}^2(\ov A,\mu_n)$ is the
direct sum of $\Ga$-modules $\NS(\ov A)/n\oplus\Br(\ov A)_n$, 
and we have
\begin{equation}
\Br(\ov A)_n=\Hom(E_n,E'_n). \label{e3}
\end{equation}

\bpr \label{1.1} Let $E$ be an elliptic curve such that the
representation of $\Ga$ in $E_\ell$ is a surjection $\Ga\to
\GL(E_\ell)$ for every prime $\ell$. Let $E'$ be an elliptic curve
with complex multiplication over $\ov k$,
which has a $k$-point of order $6$. Then
for $A=E\times E'$ we have $\Br(\ov A)^\Ga=0$. \epr {\em Proof}
Since $\Br(\ov A)$ is a torsion group it is enough to prove that for
every prime $\ell$ we have $\Br(\ov
A)_\ell^\Ga=\Hom_\Ga(E_\ell,E'_\ell)=0$.

By assumption $E'$ has complex multiplication by some imaginary
quadratic field $K$. Thus there exists an extension $k'/k$ of degree
at most 2 such that the image of ${\rm Gal}(\ov k/k')$ in
$\Aut(E'_\ell)$ is abelian. Thus the image of $\Ga$ in
$\Aut(E'_\ell)$ is a solvable group. We note that for $\ell\geq 5$
the group $\GL(2,\F_\ell)$ is not solvable. This implies that $E$
has no complex multiplication over $\ov k$. 
It follows that $E$ and $E'$ are not isogenous over $\ov k$.

The $\Ga$-module $E_\ell$ is simple, hence
any non-zero homomorphism of $\Ga$-modules $E_\ell\to E'_\ell$
must be an isomorphism. This gives a contradiction for $\ell\geq 5$.
If $\ell=2$ or $\ell=3$, the curve $E'$ has a $k$-point
of order $\ell$, so that $E'_\ell$ is not a simple $\Ga$-module,
which is again a contradiction. QED

\medskip

\noindent{\bf Example A1}
Let $k=\Q$, let $E$ be the curve $y^2=x^3+6x-2$ of
conductor $2^6 3^3$, and let $E'$ be the curve $y^2=x^3+1$ with
the point $(2,3)$ of order $6$. It follows from
\cite{Serre72}, 5.9.2, p. 318, that the conditions of Proposition
\ref{1.1} are satisfied.

\medskip

\noindent{\bf Example A2} Here is a somewhat different construction
for case A, again over $k=\Q$. Let us call a pair of elliptic curves
$(E,E')$ {\it non-exceptional} if for all primes $\ell$ the image of
the Galois group $\Ga=\Gal(\ov\Q/\Q)$ in
$\Aut(E_\ell)\times\Aut(E'_\ell)$ is as large as it can possibly be,
that is, it is the subgroup of $\GL(2,\F_\ell)\times\GL(2,\F_\ell)$
given by the condition $\det(x)=\det(x')$. This implies
$\Hom_\Ga(E_\ell,E'_\ell)=0$, so that $\Br(\ov A)^\Ga=0$, where
$A=E\times E'$. For example, let $E$ be the curve $y^2+y=x^3-x$ of
conductor 37, and let $E'$ be the curve $y^2+y=x^3+x^2$ of conductor
43. The curve $E$ has multiplicative reduction at 37, whereas $E'$
has good reduction, therefore $E$ and $E'$ are not isogenous over
$\ov\Q$. By the remark on page 329 of \cite{Serre72} the pair
$(E,E')$ is non-exceptional. In fact, most pairs $(E,E')$ are
non-exceptional in a similar sense to the remark after Proposition
\ref{pr1} (Nathan Jones, see \cite{J}).

\medskip

We now explore some other constructions providing an infinite series
of examples when $\Br(\ov A)^\Ga$ has no elements of odd order. Later we shall show that for such abelian surfaces $A$ we
often have $\Br(\Kum(A))=\Br(\Q)$, see Example 3 in Section \ref{s5}.

\bpr \label{pr} Let $E$ be an elliptic curve over $\Q$ such that
$\val_5(j(E))=-2^m$ and $\val_7(j(E))=-2^n$, where $m$ and $n$ are
non-negative integers. Let $E'$ be an elliptic curve over $\Q$ with
good reduction at $5$ and $7$, and with rational $2$-torsion, i.e.
$E'_2\subset E'(\Q)$. Then $E$ and $E'$ are not isogenous over $\ov
\Q$, and $\Hom_\Ga(E_\ell,E'_\ell)=0$ for any prime $\ell\not=2$. If
$A=E\times E'$, then $\Br(\ov A)^\Ga$ is a finite abelian $2$-group.
\epr

{\em Proof} Since $j(E)$ is not a $5$-adic integer, $E$ has
potential multiplicative reduction at $5$. But $E'$ has good
reduction at $5$, so $E$ and $E'$ are not isogenous over $\ov \Q$.

Let $p=5$. Our assumption implies that there exists a {\it Tate curve}
$\tilde E$ over $\Q_p$ such that $E\times_\Q\Q_p$ 
is the twist of $\tilde E$ by a quadratic or trivial character 
$$\chi:\Gal(\ov \Q_p/\Q_p) \to \{\pm 1\}.$$
Consider the case when $\ell\not=5$. Let $K$ be the extension
of $\Q_p$ defined as follows: if $\chi$ is trivial or unramified, then
$K=\Q_p$, and if $\chi$ is ramified, then $K\subset\ov\Q_p$ is the
invariant subfield of $\Ker(\chi)$. Let $\p$ be the maximal ideal of
the ring of integers of $K$. We note that in both cases the residue
field of $K$ is $\F_p$.

Since $\tilde{E}$ is a Tate curve, the $\ell$-torsion 
$\tilde{E}_{\ell}$ contains a Galois submodule isomorphic to 
$\mu_{\ell}$. Then the quotient $\tilde E_{\ell}/\mu_{\ell}$ is
isomorphic to the trivial Galois module $\Z/\ell$.
Hence there is a basis of $E_\ell$ such that
the image of $\Gal(\ov \Q_p/K)$ in $\Aut(E_{\ell})\simeq \GL(2,\F_\ell)$
is contained in the subgroup of upper-triangular matrices.
Let $q_E$ be the multiplicative period of $\tilde E$.
Since $\val_\p(q_E)=-\val_\p(j(E))$ is not divisible by the odd
prime $\ell$, the image of the inertia group $I(\p)$ in
$\Aut(E_{\ell})$ contains $\Id+N$ for some nilpotent $N\not=0$,
see \cite{Serre68}, Ch. IV, Section 3.2, Lemma 1. Thus
$E_{\ell}$ has exactly one non-zero $\Gal(\ov \Q_p/K)$-invariant 
subgroup $C\not=E_\ell$. As a $\Gal(\ov
\Q_p/K)$-module, $C$ is isomorphic to $\mu_{\ell}$ if $K\not=\Q_p$,
and to $\mu_\ell$ twisted by the unramified character $\chi$ if
$K=\Q_p$. The $\Gal(\ov \Q_p/K)$-module $E_{\ell}/C$ is
isomorphic to $\Z/\ell$ if $K\not=\Q_p$, and to $\Z/\ell$ twisted by
$\chi$ if $K=\Q_p$. In particular, the $\Gal(\ov \Q_p/K)$-modules 
$C$ and $E_{\ell}/C$ are isomorphic if and only if $K$ contains
a primitive $\ell$-th root of unity.

Suppose that there exists a non-zero homomorphism of $\Gal(\ov
\Q_p/K)$-modules $\phi: E_{\ell} \to E'_{\ell}$. Since $E'$ has good
reduction at $p$, the inertia $I(\p)$ acts trivially on $E'_\ell$;
in particular $\phi$ is not an isomorphism. Then $\Ker(\phi)=C$, and
$E'_{\ell}$ contains a $\Gal(\ov \Q_p/K)$-submodule isomorphic to
$\Z/\ell$ (when $\chi$ is ramified) or $\Z/\ell$ twisted by $\chi$
(when $\chi$ is unramified or trivial). 
In the first case let $E''=E'$, and in
the second case let $E''$ be the quadratic twist of $E'$ by $\chi$.
Then $E''(K)$ contains a point of order $\ell$, so that
$E''(K)$ contains a finite subgroup of order $4\ell$. Since $E'$
has good reduction at $\p$, the curve $E''$ also has good reduction
at $\p$. Then the group of $\F_5$-points on the reduction has at
least 12 elements, which contradicts the Hasse bound,
according to which an elliptic curve over $\F_p$ cannot have more
than $p+2\sqrt{p}+1$ points, see \cite{Knapp}.

It remains to consider the case $\ell=5$. The above arguments work
equally well with $p=7$. We obtain a contradiction with the Hasse
bound since no elliptic curve over $\F_7$ can contain as many as
$4\ell=20$ rational points.

The last statement of (i) follows from formula (\ref{e3}). QED

\medskip

\noindent{\bf Example A3} The curve $E$ in this proposition can
be any curve with equation $y^2=x(x-a)(x-b)$, where $a$ and $b$
are distinct non-zero integers such that exactly one of the numbers
$a$, $b$, $a-b$ is divisible by $5$, exactly one is divisible by
$7$, and none are divisible by $25$ or $49$. (For example,
$a=5+35m$, $b=7+35n$, where $m\not=2+5k$ and $n\not=4+7k$.)
The modular invariant of $E$ is
$$j(E)=2^8\frac{(a^2+b^2-ab)^3}{a^2b^2(a-b)^2},$$
which is immediate from the standard formula
$$j(E)=1728\frac{4p^3}{4p^3+27q^2}$$
for $E$ written in the Weierstrass form $y^2=x^3+px+q$, 
see \cite{Knapp}, Ch. III, Sect. 2. It follows
that $\val_5(j(E))=\val_7(j(E))=-2$. The
curve $E'$ can be any curve with equation $y^2=x(x-a')(x-b')$, where
$a'$ and $b'$ are distinct non-zero integers such that $a'$, $b'$,
$a'-b'$ are coprime to 35, e.g., $a'=35m'+1$, $b'=35n'+2$ for any
$m',\,n'\in\Z$.

\medskip

{\bf Case B.} In case B the group $\NS(\ov A)$ is freely generated
by the classes of the curves $E\times\{0\}$, $\{0\}\times E$ and the
diagonal. The image of the class of the diagonal under the map
$\NS(\ov A)\to \End(\ov E)\to \End(E_n)$ is the identity, hence
Proposition \ref{n2} gives an isomorphism of Galois modules
\begin{equation}
\Br(\ov A)_n=\End(E_n)/\Z/n, \label{e7}
\end{equation}
where $\Z/n$ is the subring of scalars in $\End(E_n)$. 

\medskip

\noindent{\bf Remark} Let $A=E\times E$, where $E$
is an elliptic curve without complex multiplication over $\ov k$,
such that the image of $\Ga$ in $\Aut(E_2)$ is $\GL(2,\F_2)$.
It is easy to check that
$\Br(\ov A)^\Ga_2=(\End(E_2)/\Z/2)^\Ga$ has order $2$; in fact,
the non-zero element of this group can be
represented by a symmetric $2\times 2$-matrix $S$ over $\F_2$
such that $S^3=\Id$. Thus the
$2$-primary component of $\Br(\ov A)^{\Gamma}$ is finite cyclic. In this
example the map $\H_{\et}^2(\ov A,\Z/2)^\Ga\to \Br(\ov A)_2^\Ga$ is
zero. By the second formula of 
Proposition \ref{n2} the map $\Br(A)_2\to
\Br(\ov A)_2^\Ga$ is {\it not} surjective.
The following proposition shows that the non-zero element
of $\Br(\ov A)_2^\Ga$ does not belong to the image of the map
$\Br(A)\to \Br(\ov A)^\Ga$.

\bpr \label{pr1}
Let $A=E\times E$, where
$E$ is an elliptic curve such that for every prime $\ell$
the image of $\Ga$ in $\Aut(E_\ell)$ is $\GL(2,\F_\ell)$.
Then we have

{\rm (i)} $\Br(A)=\Br_1(A)$;

{\rm (ii)} $\Br(\ov A)^\Ga\simeq \Z/2^m$ for some $m\geq 1$.
\epr
{\em Proof} (i) We note that the argument in the proof of
Proposition \ref{1.1} shows that the curve $E$ has no
complex multiplication. In view of the second formula of 
Proposition \ref{n2} it is enough to prove
the following lemma:

\ble \label{le2} Let $G\subset\GL(2,\Z_\ell)$ be a subgroup that
maps surjectively onto $\GL(2,\F_\ell)$. Let $\Mat_2(\Z/\ell^n)$ be
the abelian group of $2\times 2$-matrices with entries in
$\Z/\ell^n$, and let $\Mat_2(\Z/\ell^n)^G$ be the subgroup of matrices 
commuting with the image of $G$ in
$\GL(2,\Z/\ell^n)$. Then for any positive integer $n$ we have
$\Mat_2(\Z/\ell^n)^G=\Z/\ell^n\cdot\Id$. 
\ele 
{\em Proof} We proceed
by induction starting with the obvious case $n=1$. Suppose we know
the statement for $n$, and need to prove it for $n+1$. Consider the
exact sequence of $G$-modules
$$0\to\Mat_2(\Z/\ell)\to\Mat_2(\Z/\ell^{n+1})\to\Mat_2(\Z/\ell^n)\to
0,$$ where the second arrow comes from the injection
$\Z/\ell=\ell^n\Z/\ell^{n+1}\hookrightarrow \Z/\ell^{n+1}$, and the
third one is the reduction modulo $\ell^n$. By induction assumption
$\Mat_2(\Z/\ell^n)^G=\Z/\ell^n\cdot\Id$. Thus, the map
$\Mat_2(\Z/\ell^{n+1})^G\to\Mat_2(\Z/\ell^n)^G$ is surjective, and
every element in $\Mat_2(\Z/\ell^{n+1})^G$ is the sum of a scalar 
multiple of $\Id$ and an element of $\Mat_2(\Z/\ell)^G$. But
$\Mat_2(\Z/\ell)^G=\Z/\ell\cdot\Id$, and so the lemma, and 
hence also part (i) of the proposition, are proved. 

(ii) For an odd prime $\ell$ we have a direct sum decomposition of 
$\Ga$-modules $\End(E_\ell)=\Z/\ell \oplus \Br(\ov A)_\ell$, 
where $\Br(\ov A)_\ell$ is identified with the
group of endomorphisms of trace zero. Our assumption implies that 
$\Br(\ov A)^\Ga_\ell=0$. The remark before the proposition
shows that $\Br(\ov A)^\Ga$ is a finite cyclic 2-group. QED

\medskip

\noindent{\bf Remark.} By a theorem of W. Duke \cite{Duke} `almost
all' elliptic curves over $\Q$ satisfy the assumption of Proposition
\ref{pr1}. More precisely, if $y^2=x^3+ax+b$ is the unique equation
for $E$ such that $a,\,b\in\Z$ and ${\rm gcd}(a^3,b^2)$ does not
contain twelfth powers, the height $H(E)$ of $E$ is defined to be
$max(|a|^3,|b|^2)$. For $x>0$ write $\sC(x)$ for the set of elliptic
curves $E$ over $\Q$ (up to isomorphism) such that $H(E)\leq x^6$,
and $\sE(x)$ for the set of curves in $\sC(x)$ for which there
exists a prime $\ell$ such that the image of $\Ga$ in $\Aut(E_\ell)$
is not equal to $\GL(2,\F_\ell)$. Then $\lim_{x\to
+\infty}|\sE(x)|/|\sC(x)|=0$. By Proposition \ref{2.1} and Theorem
\ref{h6} this implies that for most Kummer surfaces
$$z^2=(x^3+ax+b)(y^3+ay+b)$$ we have
$\Br(X)=\Br(\Q)$. In particular, there are infinitely many such
surfaces.

\medskip

\bpr \label{pr2} Let $E$ be an elliptic curve over $\Q$ satisfying
the assumptions of Proposition $\ref{pr}$. Then $E$ has no complex
multiplication, and $\End_\Ga(E_\ell)$ is the subring of scalars
$\F_\ell\cdot {\rm Id}\subset\End(E_\ell)$ for any prime
$\ell\not=2$. If $A=E\times E$, then $\Br(\ov A)^\Ga$ is a finite
abelian $2$-group. 
\epr 
{\em Proof} The curve $E$ has no complex
multiplication because $j(E)$ is not an algebraic integer. Let
$\ell$ be an odd prime, $\ell\not=p=5$, and let $\phi$ be a non-zero
endomorphism of the $\Ga$-module $E_\ell$ such that $\Tr(\phi)=0$. From
the proof of Proposition \ref{pr} we know that there exists a nilpotent
$N\not=0$ in $\End(E_\ell)$ such that the image of $I(\p)$ in $\Aut(E_\ell)$
contains $\Id+N$.
Since $\phi$ is an endomorphism of the $I(\p)$-module $E_\ell$, it
commutes with $N$, and
it follows from $\Tr(\phi)=0$ that $\phi$ is also nilpotent.
As was explained in the proof of Proposition \ref{pr},
the existence of such an endomorphism $\phi$ implies that
$\Gal(\ov \Q_p/K)$-modules $\Z/\ell$ and $\mu_\ell$ are isomorphic.
However, $K$ does not contain non-trivial roots of 1 of order $\ell$
when $p=5$ and $\ell\not=5$ is odd, because the residue field of $K$ is $\F_5$.
This contradiction shows that $\End_\Ga(E_\ell)$ is the
subring of scalars $\Z/\ell\subset\End(E_\ell)$. If $\ell=5$ we
repeat these arguments with $p=7$ taking into account that $\Q_7$,
and hence $K$, does not contain non-trivial $5$-th roots of 1. QED

\medskip

\noindent{\bf Example B1} Let $A=E\times E$, where $E/\Q$ is an
elliptic curve such that the representation of $\Ga$ in $E_\ell$ is
a surjection $\Ga\to \GL(E_\ell)$ for every odd prime $\ell$.
Then $E$ has no complex multiplication over $\ov \Q$ (see the proof 
of Proposition \ref{1.1}).
Then $\Br(\ov A)^\Ga_\ell=0$ if $\ell\not=2$. This assumption holds in
the following examples which have the additional property that the
2-torsion of $E$ is rational, i.e., $E_2\subset E(\Q)$. For example,
one can consider the curve $y^2=(x-1)(x-2)(x+2)$ of conductor 24, or
the curve $y^2=(x+1)(x+2)(x-3)$ of conductor 40. The computation of
residues shows that in each of these two cases we have
$\Br(X)=\Br(\Q)$ for $X=\Kum(E\times E)$ (see Example 4 in the next
section).

\medskip

{\bf Case C.} Lastly, we would like to consider the case when
$A=E\times E$, where $E$ is an elliptic curve with complex
multiplication.

\bthe \label{CM} Let $E$ be an elliptic curve over $\Q$ with complex
multiplication, and let $\ell$ be an odd prime such that $E$ has no
rational isogeny of degree $\ell$, i.e., $E_{\ell}$ does not contain
a Galois-invariant subgroup of order $\ell$. Let $G_{\ell}$ be the
image of $\Gamma=\Gal(\ov\Q/\Q)$ in $\Aut(E_{\ell})$. Then
$G_{\ell}$ is nonabelian, the order $|G_\ell|$ is not divisible by
$\ell$, and the centralizer of $G_{\ell}$ in $\End(E_{\ell})$ is
$\F_{\ell}=\Z/\ell$. \ethe

The theorem remains true if one replaces $\Q$ by any real number
field (with the same proof).

\medskip

\noindent{\em Proof of Theorem} Suppose that $E$ has complex
multiplication by an order $\O$ of an imaginary quadratic field $K$,
that is, $\End(\ov E)=\O$. We start with the observation that $\ell$
is unramified in $\O$, or, equivalently, the $2$-dimensional
$\F_{\ell}$-algebra $\O/\ell \subset \End(E_{\ell})$ has no
nilpotents. Indeed, if the radical of $\O/\ell$ is non-zero, it is
an $\F_{\ell}$-vector space of dimension $1$, and so is spanned by
one element. Its kernel in $E_{\ell}$ is a Galois-invariant cyclic
subgroup of order $\ell$. We assumed that such subgroups do not
exist, so this is a contradiction.

Therefore, $\O/\ell$ is either $\F_{\ell}\oplus\F_{\ell}$ or the
field $\F_{\ell^2}$. In the first case $(\O/\ell)^{*}$ is a split
Cartan subgroup of order $(\ell-1)^2$, whereas in the second case it
is a non-split Cartan subgroup of order $\ell^2-1$, so that $\ell$
does not divide $|(\O/\ell)^{*}|$. On the other hand, the image of
$\Gal(\ov\Q/K)$ in $\Aut(E_{\ell})$ commutes with the Cartan
subgroup $(\O/\ell)^{*}$, and so belongs to $(\O/\ell)^{*}$. Since
$\Gal(\ov\Q/K)$ is a subgroup of $\Ga$ of index $2$, we conclude
that the order of $G_{\ell}$ divides $2|(\O/\ell)^{*}|$ and so is
not divisible by $\ell$.

The group $G_{\ell}$ contains an element $c$ corresponding to the
complex conjugation. Any $z\in\O\setminus\ell\O$ such that
$\Tr(z)=0$ anticommutes with complex conjugation. Since $\ell$ is
odd, the non-zero image of $z$ in $\O/\ell$ anticommutes with $c$.
Thus $c$ is not a scalar; in particular, $c$ has exact order $2$. If
$G_{\ell}$ is abelian, then both eigenspaces of $c$ in $E_{\ell}$
are Galois-invariant cyclic subgroups of order $\ell$, but these do
not exist. This implies that $G_{\ell}$ is nonabelian.

Finally, the absence of Galois-invariant order $\ell$ subgroups in
$E_{\ell}$ implies that the $G_{\ell}$-module $E_{\ell}$ is simple,
so the centralizer of $G_{\ell}$ in $\End(E_{\ell})$ is $\F_{\ell}$.
QED

\bco Let $A=E\times E$, where $E$ is an elliptic curve over $\Q$
with complex multiplication, and let $\ell$ be an odd prime such
that $E$ has no rational isogeny of degree $\ell$. Then $\Br(\ov
A)_\ell^{\Gamma}=0$. \eco {\em Proof} It follows from Theorem
\ref{CM} that the $\Ga$-module $\End(E_{\ell})$ is semisimple, hence
$\H^2_{\et}(\ov A,\mu_\ell)=(\Z/\ell)^2\oplus\End(E_\ell)$ is also
semisimple. Thus $\H^2_{\et}(\ov A,\mu_\ell)=\NS(\ov
A)/\ell\oplus\Br(\ov A)_\ell$ is a direct sum of $\Ga$-modules.
Since the identity in $\End(E_\ell)$ corresponds to the diagonal in
$\ov E\times\ov E$, it is contained in $\NS(\ov A)/\ell$. By Theorem
\ref{CM} we have $\H^2_{\et}(\ov A,\mu_\ell)^\Ga\subset\NS(\ov
A)/\ell$, so that $\Br(\ov A)_\ell^{\Gamma}=0$. QED

\medskip

\noindent{\bf Example C1} Let $A=E\times E$, where $E$ is the curve
$y^2=x^3-x$ with complex multiplication by $\Z[\sqrt{-1}]$. An
application of {\tt sage} \cite{sage} gives that every isogeny
of prime degree $E\to E'$ defined over $\Q$
is the factorization by a subgroup of
$E(\Q)_{\rm tors}=E_2$. Hence $\Br(\ov A)_\ell^{\Gamma}=0$ for every
odd prime $\ell$.

\medskip

\noindent{\bf Example C2} Let $A=E\times E$, where $E$ is the curve
$y^2=x^3-1$ with complex multiplication by
$\Z[\frac{1+\sqrt{-3}}{2}]$. An application of {\tt sage} gives that
every isogeny of prime degree $E\to E'$ over $\Q$
is the factorization by a subgroup
of $E(\Q)_{\rm tors}\simeq\Z/6$. Hence $\Br(\ov A)_\ell^{\Gamma}=0$
for every prime $\ell\geq 5$.

\section{Brauer groups of Kummer surfaces} \label{s5}

\noindent{\bf Example 1} Let $k=\Q$. Examples A1 and A2 show that
the Kummer surfaces $X$ given by the following affine equations have
trivial Brauer group $\Br(X)=\Br(\Q)$:
\begin{equation}
z^2=(x^3+6x-2)(y^3+1),\label{ex}
\end{equation}
\begin{equation}
z^2=(4x^3-4x+1)(4y^3+4y^2+1).\label{ex2}
\end{equation}
In both examples we have $\Br(\ov X)^\Ga=0$.

\noindent{\bf Example 2} Other examples can be obtained using
Proposition \ref{pr1} in conjunction with Theorem \ref{h6}. For
example, for the following Kummer surface $X$ we also have
$\Br(X)=\Br(\Q)$, whereas $\Br(\ov X)^\Ga\simeq\Z/2^m$ for some $m\geq 1$:
\begin{equation}
z^2=(x^3+6x-2)(y^3+6y-2).\label{ex3}
\end{equation}

The interest of the following series of examples is that for them
the image of $\Br(A)$ in $\Br(\ov A)$ contains $\Br(\ov A)_2$, so in
order to prove the triviality of $\Br(X)$ we need to compute the
residues at the nine lines in $X\setminus W$.

\smallskip

\noindent{\bf Example 3} Let $X$ be the Kummer surface over $\Q$
with affine equation
$$z^2=x(x-a)(x-b)y(y-a')(y-b'),$$ such that $a=5+35m$, $b=7+35n$, where
$m,\,n\in \Z$, $m$ is not congruent to $2$ modulo $5$, $n$ is not
congruent to $4$ modulo $7$, and $a'=35m'+1$, $b'=35n'+2$ for any
$m',\,n'\in\Z$. We have $X=\Kum(E\times E')$, where the elliptic
curves $E$ and $E'$ are as in Example A3. 
Since $X(\Q)\not=\emptyset$ we see that $\Br(\Q)$
is a direct factor of $\Br(X)$.
By Propositions \ref{pr}
and \ref{2.1} (ii) to show that  $\Br(X)=\Br(\Q)$ it is enough to
prove that every element of $\Br(X)$ of order 2 is algebraic. By
Proposition \ref{res1} we need to compute the dimension $d$ 
of the kernel 
of the matrix (\ref{res}). Considering the first two entries in each
row, and taking their valuations at 5 and 7 immediately shows that
no product of some of the rows of (\ref{res}) is trivial. Thus 
$d=0$, hence $\Br(X)=\Br(\Q)$.

\medskip

\noindent{\bf Example 4} Let $X=\Kum(E\times E)$, where $E$ is as in
Example 3, or the elliptic curve with conductor 24 or 40 mentioned
in Example B1. In the latter case $X$ is given by one of the
following equations:
\begin{equation}
z^2=(x-1)(x-2)(x+2)(y-1)(y-2)(y+2),\label{ex4}
\end{equation}
\begin{equation}
z^2=(x+1)(x+2)(x-3)(y+1)(y+2)(y-3).\label{ex5}
\end{equation}
One checks that
the dimension of the kernel of (\ref{res}) is 1, so that $\Br(X)=\Br(\Q)$ by
Proposition \ref{res1}.

\bigskip

\noindent{\bf Kummer surfaces without rational points}. There is a
more general construction of Kummer surfaces than the one previously
considered. Let $c$ be a 1-cocycle of $\Ga$ with coefficients in
$A_2$ so that $[c]\in\H^1(k,A_2)$. All quasi-projective varieties and Galois
modules acted on by the $k$-group scheme
$A_2$ can be twisted by $c$. The twist of $A$
is a principal homogeneous space $A^c$, also called a 2-covering
of $A$. The action of $A_2$ on $A$ by translations descends to an
action of $A_2$ on $X=\Kum(A)$, so we obtain a twisted Kummer
surface $X^c$ together with morphisms 
$A^c\leftarrow {A'}^c\to X^c$. For example, if $A=E\times E'$, 
a 2-covering $C$ of $E$
is given by $y^2=f(x)$, where $f(x)$ is a separable polynomial of
degree 4, and a 2-covering $C'$ of $E'$ is given by a similar equation
$y^2=g(x)$, then the
twisted Kummer surface $X^c$ is given by the affine equation
$$ z^2=f(x)g(y).$$ The Hasse
principle on such surfaces over number fields was studied in
\cite{SS}.

\bpr \label{t1}
Suppose that for every integer $n>1$ we have 
$\H^2_{\et}(\ov A,\mu_n)^\Ga=(\NS(\ov A)/n)^\Ga$
(this condition is satisfied when $\Br(\ov A)^\Ga=0$).
Then $\Br(X^c)=\Br_1(X^c)$ for
any $[c]\in\H^1(k,A_2)$.
\epr
{\em Proof} By Remark 3 after
Proposition \ref{2.5}, $A_2(\ov k)$ acts trivially on $\Br(\ov A)$
and on $\Br(\ov X)$,
thus we have the following isomorphisms of $\Ga$-modules:
$$\Br(\ov{X^c})\simeq \Br(\ov X)\simeq
\Br(\ov{A^c})\simeq\Br(\ov A).$$ Translations act trivially
on \'etale cohomology groups of $\ov A$, hence we have a canonical
isomorphism of $\Ga$-modules
$\H_{\et}^2(\ov {A^c},\mu_{n})=\H_{\et}^2(\ov {A},\mu_{n})$.
In the commutative diagram
$$\begin{array}{ccc}
\H_{\et}^2(\ov {A^c},\mu_{n})^\Ga&\to &\Br(\ov {A^c})_{n}^\Ga\\
\uparrow&&\uparrow\\
\H_{\et}^2(A^c,\mu_{n})&\to &\Br(A^c)_{n}
\end{array}$$
the bottom arrow is surjective, and the top arrow is zero by
assumption. It follows that $\Br(A^c)=\Br_1(A^c)$. 
We conclude by Theorem \ref{h6}. QED

\medskip

This proposition in conjunction with Proposition \ref{1.1} 
gives many examples of twisted Kummer surfaces $X^c$
such that $\Br(X^c)=\Br_1(X^c)$.

\bigskip

\noindent{\bf Kummer surfaces with non-trivial transcendental Brauer
group}. Let $E$ be an elliptic curve over $\Q$. As pointed out by
Mazur \cite[p. 133]{Mazur} the elliptic curves $E'$ such that the
Galois modules $E_\ell$ and $E'_\ell$ are symplectically isomorphic
correspond to $\Q$-points on the modular curve $X(\ell)$ twisted by
$E_\ell$. Thus for $\ell\leq 5$ there are infinitely many
possibilities for $E'$ due to the fact that the genus of $X(\ell)$
is zero, see \cite{Si}, \cite{RS}. Now let $\ell=7$, $11$ or $13$.
Examples of pairs of non-isogenous elliptic curves with
isomorphic Galois modules $E_\ell\simeq E'_\ell$ for these values
of $\ell$ can be found in \cite{H} and \cite{D}. 
Our results imply that $\Br(X)$, where $X=\Kum(E\times E')$, contains an
element of order $\ell$ whose image in $\Br(\ov X)$ is non-zero.

\bigskip

\noindent{\it Acknowledgements.} The first named author (A.S.) thanks
Nicolas Billerey, Nathan Jones and Christian Wuthrich
for useful discussions and their help with references. 
The use of {\tt sage} is gratefully acknowledged. 
The work on this paper started at the
workshop `The arithmetic of K3 surfaces' at the Banff International
Research Station, attended by both authors,
and was mostly done while A.S.
enjoyed the hospitality of the Hausdorff Research Institute for
Mathematics in Bonn. A.S. thanks the organizers of the
special semester on Diophantine equations for support and excellent
working conditions. The authors are deeply grateful to Evis Ieronymou
for his careful reading of the paper and helpful remarks.

\bigskip

\noindent Department of Mathematics, South Kensington Campus,
Imperial College London, SW7 2BZ England, U.K.

\smallskip

\noindent Institute for the Information Transmission Problems,
Russian Academy of Sciences, 19 Bolshoi Karetnyi, Moscow, 127994
Russia
\medskip

\noindent a.skorobogatov@imperial.ac.uk

\bigskip

\noindent Department of Mathematics, Pennsylvania State University,
University Park, Pennsylvania 16802, USA

\smallskip

\noindent Institute for Mathematical Problems in Biology, Russian
Academy of Sciences, Pushchino, Moscow Region, Russia

\noindent zarhin@math.psu.edu

\end{document}